\documentclass[12pt]{amsart} 
\usepackage{amsmath,amssymb}
\usepackage{latexsym}
\usepackage{amssymb} % \mathbb
\usepackage{latexsym}
\usepackage{delarray} 

\def\R{{\mathbb R}}

\def\Rn{{{\mathbb R}^n}}

\def\Rnx{{\R^n_x}}

\def\L2tx{{L^2(\R_t\times\R^n_x)}}
\def\Lx{{L^2(\R^n_x)}}
\def\p#1{{\left({#1}\right)}}

\def\n#1{{\left\|{#1}\right\|}}
\def\abs#1{{\left|{#1}\right|}}
\def\jp#1{{\left\langle{#1}\right\rangle}}

\def\supp{\operatorname{supp}}

\def\va{\varphi}

\numberwithin{equation}{section}
\theoremstyle{plain}
\newtheorem{thm}{Theorem}[section]

\newtheorem{cor}[thm]{Corollary}

\begin{document}

\title{Comparison of estimates
for dispersive equations}

\author{M. Ruzhansky}

\address{Department of Mathematics, Imperial College London\\
London, United Kingdom\\
E-mail: m.ruzhansky@imperial.ac.uk}

\author{M. Sugimoto}

\address{Department of Mathematics, Osaka University\\
Osaka, Japan\linebreak
E-mail: sugimoto@math.sci.osaka-u.ac.jp}

\begin{abstract}
This paper describes a new comparison principle that
can be used for the comparison
of space-time estimates for dispersive equations.
In particular, results are applied to the global smoothing
estimates for several classes of 
dispersive partial differential equations.
\end{abstract}


\maketitle

%\bodymatter

\section{Introduction}\label{aba:sec1}

In this note we will present a new comparison principle that
allows one to compare certain estimates for dispersive equations
of different types based on expressions involving their
symbols. In particular, a question is that if we have a
certain estimate for one equation, whether we can derive
a corresponding estimate for another equation. This question
is of interest on its own, and it has several
applications. Proofs of the statements of this paper 
can be found in the authors' paper \cite{RSarx}.

The main 
application of this technique that we have in mind is for
the global smoothing estimates for dispersive equations. These
smoothing estimates are essentially global space-time
estimates in weighted Sobolev spaces over $L^2$,
see for example
\cite{BD91, BK92, CS88, Ka83, Sj87, Ve88}. There
is a known method in the microlocal analysis on how to
transform one equation into another, namely the canonical
transforms realised in the form of Fourier integral
operators \cite{DH72}. In the global setting one needs to
develop global weighted estimates in $L^2$ for the 
corresponding classes of Fourier integral operators
in order to apply them to the smoothing estimates. Such global
estimates for Fourier integral operators have been established
by the authors \cite{RS06a} and have been applied to derive
new smoothing estimates for Schr\"odinger equations
\cite{RS06a, RS06}. These techniques allow one to reduce the
analysis of dispersive equations to normal forms in one
and two dimensions \cite{RS04}. The comparison principles
introduced below allow one to further relate estimates in
normal forms, thus establishing comprehensive relations
between smoothing estimates for dispersive equations with
constant coefficient \cite{RSarx}.

In this note, we denote 
$x=(x_1,\ldots,x_n)$, $\xi=(\xi_1,\ldots,\xi_n)$, and
$D_x=(D_1,D_2\ldots,D_n)$, where
$D_j$ denotes $D_{x_j}=\frac{1}{i}\frac{\partial}
{\partial x_j}$, for all $j=1,2,\ldots,n$, and
$i=\sqrt{-1}$.

\section{Comparison principles}

Let $u(t,x)=e^{it f(D_x)}\varphi(x)$ and 
$v(t,x)=e^{it g(D_x)}\varphi(x)$
be solutions to the following evolution equations, 
where $t\in\R$ and $x\in\R^n$:
\begin{equation}\label{EQ:maineq1}
\left\{
\begin{aligned}
\p{i\partial_t+f(D_x)}\,u(t,x)&=0,\\
u(0,x)&=\varphi(x),
\end{aligned}
\right.
\end{equation}
and 
\begin{equation}\label{EQ:maineq2}
\left\{
\begin{aligned}
\p{i\partial_t+g(D_x)}\,v(t,x)&=0,\\
v(0,x)&=\varphi(x).
\end{aligned}
\right.
\end{equation}

First we state the following result relating several norms
involving propagators for equations \eqref{EQ:maineq1}
and \eqref{EQ:maineq2}:
\medskip
\begin{thm}\label{prop:basiceq}
Let $f\in C^1(\R^n)$ be a real-valued function
such that, for almost all $\xi'=(\xi_2,\ldots,\xi_n)\in\R^{n-1}$,
$f(\xi_1,\xi')$ is strictly monotone in $\xi_1$ 
on the support of a measurable function $\sigma$ on $\R^n$.
Then we have
\begin{equation*}\label{EQ:basiceq}
\n{\sigma(D_x)
e^{it f(D_x)}\varphi(x_1,x')}_{L^2(\R_t\times\R_{x'}^{n-1})}^2
= (2\pi)^{-n}
 \int_\Rn |\widehat{\varphi}(\xi)|^2 \frac{|\sigma(\xi)|^2}
 {|\partial f/\partial\xi_1(\xi)|}\, d\xi
\end{equation*}
for all $x_1\in\R$, where $x'=(x_2,\ldots,x_n)\in\R^{n-1}$.
\end{thm}

The following comparison principle
is a straightforward consequence of Theorem \ref{prop:basiceq}:
\medskip
\begin{cor}\label{prop:dimneq}
Let $f,g\in C^1(\R^n)$ be real-valued functions
such that, for almost all $\xi'=(\xi_2,\ldots,\xi_n)\in\R^{n-1}$,
$f(\xi_1,\xi')$ and $g(\xi_1,\xi')$
are strictly monotone in $\xi_1$ 
on the support of a measurable function $\chi$ on $\R^n$.
Let $\sigma,\tau\in C^0(\R^n)$ be such that, for some $A>0$, we have
\begin{equation}\label{EQ:compassdimn}
\frac{|\sigma(\xi)|}{\abs{\partial_{\xi_1} f(\xi)}^{1/2}}
\leq A \, \frac{|\tau(\xi)|}{\abs{\partial_{\xi_1} g(\xi)}^{1/2}}
\end{equation} 
for all $\xi\in\supp\chi$ satisfying
$D_1 f(\xi)\not=0$ and 
$D_1 g(\xi)\not=0$.
Then we have
\begin{multline}\label{EQ:dimnest}
\n{\chi(D_x)\sigma(D_x)
e^{it f(D_x)}\varphi(x_1,x')}_{L^2(\R_t\times\R_{x'}^{n-1})} 
\\ \leq  A 
\|\chi(D_x)\tau(D_x)e^{it g(D_x)}\varphi(\widetilde x_1,x')\|_
{L^2(\R_t\times\R_{x'}^{n-1})}
\end{multline}
for all $x_1,\widetilde x_1\in\R$, where $x'=(x_2,\ldots,x_n)\in\R^{n-1}$.
Consequently, for any  measurable function $w$ on $\R$ we have
\begin{multline}\label{EQ:dimnestwgt}
\n{w(x_1)\chi(D_x)\sigma(D_x)
e^{it f(D_x)}\varphi(x)}_{L^2(\R_t\times\R_{x}^{n})} 
\\ \leq  A 
\|w(x_1)\chi(D_x)\tau(D_x)e^{it g(D_x)}\varphi(x)\|_
{L^2(\R_t\times\R_{x}^{n})}.
\end{multline}
Moreover, if $\chi\in C^0(\Rn)$ and $w\neq0$ on a set
of $\R$ with positive measure, the converse is true,
namely, if we have estimate \eqref{EQ:dimnest} for all $\va$,
for some $x_1, \widetilde{x}_1\in\R$,
or if we have estimate \eqref{EQ:dimnestwgt} for all $\va$,
and the norms are finite, then we also have
inequality \eqref{EQ:compassdimn}.
\end{cor}
\medskip
We remark that inequality \eqref{EQ:dimnestwgt}
in Corollary \ref{prop:dimneq}
gives the comparison between
different weighted estimates.
The reason to introduce
a cut-off function $\chi$ into the estimates is that the relation
between symbols may be different for different regions
of the frequencies $\xi$
(for example this is the case for the relativistic Schr\"odinger
and for the Klein-Gordon equations, which are of order two for 
large frequencies and of order zero for small
frequencies). In such case we can use this comparison principle
to relate estimates for the corresponding ranges of
frequencies, thus yielding
more refined results, since then we have freedom to
choose different $\sigma$ for different types of behaviour
of $f^\prime$.
The assumption $\sigma,\tau\in C^0(\R^n)$ that was made 
in Corollary \ref{prop:dimneq} is for
the clarity of the exposition and can clearly be
relaxed. 
\par
In the case $n=1$, we neglect $x'=(x_2,\ldots,x_n)$ 
in a natural way and
just write $x=x_1$, $\xi=\xi_1$, and $D_x=D_1$.
Similarly in the case $n=2$, we use the notation $(x,y)=(x_1,x_2)$,
$(\xi,\eta)=(\xi_1,\xi_2)$, and $(D_x,D_y)=(D_1,D_2)$.
In both cases, we write $\widetilde x=\widetilde x_1$ 
in the notation of
Corollary \ref{prop:dimneq}.
Then we have the following corollaries in lower dimensions:
\medskip
\begin{cor}\label{prop:dim1eq}
Suppose $n=1$.
Let $f,g\in C^1(\R)$ be real-valued and strictly
monotone on the support of a measurable function $\chi$ on $\R$.
Let $\sigma,\tau\in C^0(\R)$ be such that, for some $A>0$, we have
\begin{equation}\label{dim1eq:1}
\frac{|\sigma(\xi)|}{|f^\prime(\xi)|^{1/2}}
\leq A \frac{|\tau(\xi)|}{|g^\prime(\xi)|^{1/2}}
\end{equation}
for all $\xi\in\supp\chi$ satisfying
$f^\prime(\xi)\not=0$ and $g^\prime(\xi)\not=0$.
Then we have
\begin{equation}\label{dim1eq:2}
\|\chi(D_x)\sigma(D_x)e^{it f(D_x)}\varphi(x)\|_{L^2(\R_t)}\leq A
\|\chi(D_x)\tau(D_x)e^{it g(D_x)}\varphi(\widetilde x)\|_{L^2(\R_t)}
\end{equation}
for all $x,\widetilde x\in\R$.
Consequently, for general $n\geq1$ and
for any measurable function $w$ on $\Rn$, we have
\begin{multline}\label{dim1eq:3}
\|w(x)\chi(D_j)\sigma(D_j)e^{it f(D_j)}\varphi(x)\|_
{L^2(\R_t\times\R_x^n)} \\ \leq A
\|w(x)\chi(D_j)\tau(D_j)e^{it g(D_j)}\varphi(x)\|_
{L^2(\R_t\times\R_x^n)},
\end{multline}
where $j=1,2,\ldots,n$.
Moreover, if $\chi\in C^0(\R)$ and $w\neq0$ on a set
of $\R^n$ with positive measure,
the converse is true,
namely, if we have estimate \eqref{dim1eq:2} for all $\va$,
for some $x, \widetilde{x}\in\R$,
or if we have estimate \eqref{dim1eq:2} for all $\va$,
and the norms are finite, then we also have
inequality \eqref{dim1eq:1}.
\end{cor}
We have the following comparison principle in two dimensions:
\begin{cor}\label{prop:dim2eq}
Suppose $n=2$.
Let $f,g\in C^1(\R^2)$ be real-valued functions
such that, for almost all $\eta\in\R$,
$f(\xi,\eta)$ and $g(\xi,\eta)$ are strictly
monotone in $\xi$ on the support of a measurable function $\chi$ on $\R^2$.
Let $\sigma,\tau\in C^0(\R^2)$ be such that, for some $A>0$, we have
\begin{equation}\label{dim2eq:1}
\frac{|\sigma(\xi,\eta)|}
{\abs{\partial f/\partial \xi(\xi,\eta)}^{1/2}}
\leq A \frac{|\tau(\xi,\eta)|}
{\abs{\partial g/\partial \xi(\xi,\eta)}^{1/2}}
\end{equation}
for all $(\xi,\eta)\in\supp\chi$ satisfying
$\partial f/\partial \xi(\xi,\eta)\not=0$ and 
$\partial g/\partial \xi(\xi,\eta)\not=0$.
Then we have
\begin{multline}\label{dim2eq:2}
\n{\chi(D_x,D_y)\sigma(D_x,D_y)
e^{it f(D_x,D_y)}\varphi(x,y)}_{L^2(\R_t\times\R_{y})} \\ \leq A
\|\chi(D_x,D_y)\tau(D_x,D_y)e^{it g(D_x,D_y)}\varphi(\widetilde x,y)\|_
{L^2(\R_t\times\R_{y})}
\end{multline}
for all $x,\widetilde x\in\R$.
Consequently, for general $n\geq2$ and 
for any  measurable function $w$ on $\R^{n-1}$ we have
\begin{multline}\label{dim2eq:3}
\|w(\check{x}_k)\chi(D_j,D_k)\sigma(D_j,D_k)
e^{it f(D_j,D_k)}\varphi(x)\|_
{L^2(\R_t\times\R_x^n)} \\ \leq A
\|w(\check{x}_k)\chi(D_j,D_k)\tau(D_j,D_k)
e^{it g(D_j,D_k)}\varphi(x)\|_
{L^2(\R_t\times\R_x^n)},
\end{multline}
where $j\neq k$ and $\check{x}_k=(x_1,\ldots,x_{k-1},x_{k+1},\ldots,x_n)$.
Moreover, if $\chi\in C^0(\R^2)$ and $w\neq0$ on a set
of $\R^{n-1}$ with positive measure, the converse is true,
namely, if we have estimate \eqref{dim2eq:2} for all $\va$,
for some $x, \widetilde{x}\in\R$,
or if we have estimate \eqref{dim2eq:2} for all $\va$,
and the norms are finite, then we also have
inequality \eqref{dim2eq:1}.
\end{cor}
%\medskip
By the same argument as used in the proof of Theorem \ref{prop:basiceq}
and Corollary \ref{prop:dimneq},
we have a comparison result in the radially symmetric case.
Below, we denote the set of the positive real numbers $(0,\infty)$
by $\R_+$.
%\medskip
\begin{thm}\label{prop:dim1eqmod}
Let $f,g\in C^1(\R_+)$ be real-valued and strictly
monotone on the support of a measurable function $\chi$ on $\R_+$.
Let $\sigma,\tau\in C^0(\R_+)$ be such that, for some $A>0$, we have
\begin{equation}\label{EQ:compassdimmod}
\frac{|\sigma(\rho)|}{|f^\prime(\rho)|^{1/2}}
\leq A \frac{|\tau(\rho)|}{|g^\prime(\rho)|^{1/2}}
\end{equation}
for all $\rho\in\supp\chi$ satisfying
$f^\prime(\rho)\not=0$ and $g^\prime(\rho)\not=0$.
Then we have
\begin{equation}\label{EQ:dim1estmod}
\|\chi(|D_x|)\sigma(|D_x|)e^{it f(|D_x|)}\varphi(x)\|_{L^2(\R_t)}
\leq A
\|\chi(|D_x|)\tau(|D_x|)e^{it g(|D_x|)}\varphi(x)\|_{L^2(\R_t)}
\end{equation}
for all $x\in\R^n$.
Consequently, for any measurable function $w$ on $\Rn$, we have
\begin{multline}\label{EQ:dim1estxcor}
\|w(x)\chi(|D_x|)\sigma(|D_x|)e^{it f(|D_x|)}\varphi(x)\|_
{L^2(\R_t\times\R_x^n)} \\ \leq A
\|w(x)\chi(|D_x|)\tau(|D_x|)e^{it g(|D_x|)}\varphi(x)\|_
{L^2(\R_t\times\R_x^n)}.
\end{multline}
Moreover, if $\chi\in C^0(\R_+)$ and $w\neq0$ on a set
of $\R^n$ with positive measure, the converse is true,
namely, if we have estimate \eqref{EQ:dim1estmod} for all $\va$,
for some $x\in\Rn$,
or if we have estimate \eqref{EQ:dim1estxcor} for all $\va$,
and the norms are finite, then we also have
inequality \eqref{EQ:compassdimmod}.
\end{thm}
Theorem \ref{prop:dim1eqmod} provides an analytic alternative
to computations for certain estimates in the radially symmetric
case done with the help of special functions \cite{Wa02}.

These comparison principles can be extended to provide the
relation between Strichartz type norms, and the details
and the meaning of the corresponding estimates
can be found in authors' paper \cite{RSarx}. Here we just
give one corollary: 

\begin{cor}\label{cor:Strichartz}
Let functions $f,g,\sigma,\tau$ be as in Theorem 
{\rm \ref{prop:dim1eqmod}} and satisfy relation
{\rm \eqref{EQ:compassdimmod}}. Let $0<p\leq\infty$. Then, 
for any measurable function $w$ on $\Rn$, we have the
estimate
\begin{multline}\label{EQ:dim1estxcor-Str}
\|w(x)\chi(|D_x|)\sigma(|D_x|)e^{it f(|D_x|)}\varphi(x)\|_
{L^p(\Rnx,L^2(\R_t))} \\ \leq A
\|w(x)\chi(|D_x|)\tau(|D_x|)e^{it g(|D_x|)}\varphi(x)\|_
{L^p(\Rnx,L^2(\R_t))}.
\end{multline}
\end{cor}
From this, it follows, for example, that 
for all $0<p\leq\infty$,
quantities
$||e^{it\sqrt{-\Delta}}\va||_{L^p(\Rnx,L^2(\R_t))}$,
$|||D_x|^{1/2}e^{-it\Delta}\va||_{L^p(\Rnx,L^2(\R_t))}$,
and 
$|||D_x|e^{it(-\Delta)^{3/2}}\va||_{L^p(\Rnx,L^2(\R_t))}$
for the propagators of the wave, Schr\"odinger, and KdV 
type equations
are equivalent.

%\newpage
\section{Some applications}
\label{SECTION:model}
Let us now give some examples of the use of these comparison
principles. 
If both sides in expression \eqref{EQ:compassdimn} 
in Corollary \ref{prop:dimneq}
are equivalent,
we can use the comparison in two directions, from which it
follows that norms on both sides in \eqref{EQ:dimnest}
are equivalent.
The same is true for Corollaries \ref{prop:dim1eq}, \ref{prop:dim2eq}
and Theorem \ref{prop:dim1eqmod}.
In particular, we can
conclude that many smoothing estimates for the Schr\"odinger
type equations of different orders are equivalent to each other.
Indeed, applying Corollary \ref{prop:dim1eq} in two directions, we
immediately obtain that for 
$n=1$ and $l,m>0$, we have
\begin{equation}\label{prop:dim1ex}
\n{|D_x|^{(m-1)/2}e^{it|D_x|^{m}}
\varphi(x)}_{L^2(\R_t)}=
\sqrt{\frac{l}{m}}
\n{|D_x|^{(l-1)/2}e^{it|D_x|^{l}}
\varphi(x)}_{L^2(\R_t)}
\end{equation}
for every $x\in\R$,
assuming that $\supp\widehat{\varphi}\subset [0,+\infty)$ or
$(-\infty,0]$.
On the other hand, still in the case $n=1$, we have easily
\begin{equation}\label{core}
\n{e^{itD_x}\varphi(x)}_{L^2(\R_t)}
=\n{\varphi}_{L^2\p{\R_{x}}}
\quad \text{for all $x\in\R$},
\end{equation}
which is a straightforward consequence of the fact
$e^{itD_x}\varphi(x)=\varphi(x+t)$.
These observations yield:
\begin{thm}\label{prop:basic}
Suppose $n=1$ and $m>0$.
Then we have
\begin{equation}\label{basic:1dim}
\n{|D_x|^{(m-1)/2}e^{it|D_x|^m}\varphi(x)}_{L^2(\R_t)}
\leq C\n{\varphi}_{L^2(\R_x)}
\end{equation}
for all $x\in\R$.
Suppose $n=2$ and $m>0$.
Then we have
\begin{equation}\label{basic:2dim}
\n{|D_y|^{(m-1)/2}e^{itD_x|D_y|^{m-1}}\varphi(x,y)}_{L^2(\R_t\times\R_y)}
\leq C\n{\varphi}_{L^2\p{\R^2_{x,y}}}
\end{equation}
for all $x\in\R$.
Each estimate above is equivalent to itself with $m=1$
which is a direct consequence of equality \eqref{core}.
\end{thm}
\medskip
Estimates \eqref{basic:1dim} and \eqref{basic:2dim}
in Theorem \ref{prop:basic} in the special case $m=2$
were shown by Kenig, Ponce and Vega \cite[p.56]{KPV91} and by
Linares and Ponce \cite[p.528]{LP93}, respectively.
Theorem \ref{prop:basic} shows that
these results, together with their generalisation to
other orders $m$,
are in fact just corollaries of the elementary one 
dimensional fact
$e^{itD_x}\varphi(x)=\varphi(x+t)$ once we apply the comparison principle.
\par
By using the comparison principle in the radially symmetric 
and higher dimensional cases,
we have also another type of equivalence of
smoothing estimates, which can be found in authors' paper
\cite{RSarx}. Let us give one example:
\begin{thm}\label{prop:basic2}
For $m>0$ {\rm(}and any $\alpha,\beta${\rm)}
we have the following relations {\rm(}in the first equality
the left and the right hand sides are finite
for the same values of $\alpha,\beta$ at the same time{\rm)}
\begin{align*}
&\n{|x|^{\beta-1} |D_x|^{\beta} e^{it |D_x|^2}\varphi}_\L2tx
 = \\ &
\qquad\qquad\qquad\qquad\qquad\qquad \sqrt{\frac{m}{2}}
\n{|x|^{\beta-1} |D_x|^{m/2+\beta-1} e^{it |D_x|^m}\varphi}_\L2tx,
\\
&\begin{aligned}
\n{
 \jp{x}^{\alpha-m/2}|D_x|^\alpha
  e^{it|D_x|^m}\varphi(x)
}_{L^2(\R_t\times\R^n_x)}
& \\ \qquad\qquad\qquad\qquad \leq
\n{
 |x|^{\alpha-m/2}|D_x|^\alpha
  e^{it|D_x|^m}\varphi(x)
}_{L^2(\R_t\times\R^n_x)} \\
  \qquad\qquad \leq 
\sup_{\lambda>0}\n{
 \jp{x}^{\alpha-m/2}|D_x|^\alpha
  e^{it|D_x|^m}\varphi_\lambda(x)
}_{L^2(\R_t\times\R^n_x)},
\end{aligned}
\end{align*}
where $\varphi_\lambda(x)=\lambda^{n/2}\varphi(\lambda x)$, and
we take $\alpha\leq m/2$ in the last estimate.
The operator norms of operators
$\jp{x}^{\alpha-m/2} |D_x|^\alpha e^{it|D_x|^m}$ and
$|x|^{\alpha-m/2} |D_x|^\alpha e^{it|D_x|^m}$ as mappings from
$L^2(\Rn)$ to $\L2tx$ are equal.
\end{thm}
As a nice consequence, for $n\geq 3$ and $m>0$
we can conclude also the estimate
\begin{equation}\label{EQ:Simon-const}
\n{|x|^{-1} |D_x|^{m/2-1}
  e^{it|D_x|^m}\varphi(x)
}_{L^2(\R_t\times\R^n_x)}\leq
\sqrt{\frac{2\pi}{m(n-2)}}\n{\va}_{L^2(\R_x^n)},
\end{equation}
where the constant $\sqrt{\frac{2\pi}{m(n-2)}}$ is sharp.
This follows from the first equality in Theorem \ref{prop:basic2}
with $\beta=0$ and the best constant in the 
case $m=2$ (as shown by Simon \cite{Si92} as a 
consequence of constants in Kato's theory).

As a consequence of Theorem \ref{prop:basic2}, we have
\begin{cor}\label{Th:typeII}
Suppose $m>0$ and $(m-n)/2<\alpha<(m-1)/2$.
Then we have
\begin{equation}\label{model:5}
\n{\abs{x}^{\alpha-m/2}|D_x|^{\alpha}e^{it|D_x|^{m}}
\varphi(x)}_
{L^2\p{\R_t\times\R^n_x}}
\leq C\n{\varphi}_{L^2\p{\R^n_x}}.
\end{equation}
Suppose $m>0$ and $(m-n+1)/2<\alpha<(m-1)/2$.
Then we have
\begin{equation}\label{model:6}
\n{\abs{x}^{\alpha-m/2}|D'|^{\alpha} 
  e^{it\p{|D_1|^m-|D'|^m}}\varphi(x)
}_\L2tx\leq C\n{\varphi}_\Lx,
\end{equation}
where $D'=(D_2,\ldots,D_n)$.
\end{cor}
Estimate \eqref{model:5} is known in the case $m=2$ as the
Kato--Yajima estimate \cite{KY89}. The application of the
comparison principles also yields some refinements for
other equations, for example for the relativistic
Schr\"odinger equation \cite{BN97}, Klein-Gordon equations
and wave equations \cite{Be94}. We refer to \cite{RSarx}
for further details.

\bibliographystyle{ws-procs9x6}
%\bibliography{ws-pro-sample}

\end{document}